\documentclass[a4paper,11pt]{article}
\usepackage{latexsym}
\usepackage[namelimits,intlimits,sumlimits]{amsmath}
\usepackage{amsfonts}
\usepackage{amssymb}
\usepackage{epsfig}

\newtheorem{satz}{Theorem}
\newtheorem{lemma}{Lemma}
\newtheorem{rem}{Remark}
\newtheorem{kor}{Corollary}
\newtheorem{prop}{Proposition}
\newtheorem{defin}{Definition}

\newcommand{\R}{\mathbb{R}}

\newcommand{\Q}{\mathbb{Q}}
\newcommand{\Z}{\mathbb{Z}}
\newcommand{\N}{\mathbb{N}}
\newcommand{\T}{\mathbb{T}}
\renewcommand{\H}{\mathbb{H}}
\newcommand{\V}{\mathbb{V}}
\newcommand{\CC}{\mathcal{C}}

\DeclareMathOperator{\ord}{ord}

\newenvironment{proof}{{\it Proof}: }{ $\mbox{ }$ \\ $\mbox{ }$
\hfill $\square$ \newline}
\newenvironment{beweis}{{\it Proof}: }{ $\mbox{ }$ \\ $\mbox{ }$
\hfill $\square$ \newline}

\title{ On the asymptotic quadratic growth rate of saddle connections and 
periodic orbits on marked flat tori}

\author{Martin Schmoll}

\begin{document}

\maketitle

\begin{abstract}
Asymptotic quadratic growth rates of saddle connections and families 
of periodic cylinders on translation tori with $n$ marked points are studied. 
For any marking the existence of limits of the quadratic growth rate is shown 
using elementary methods (not Ratners theorem). 
We study the growth rate limit as function of the marking.
Precise formulas for this function in the case of two marked points are given and 
the sets where the growth rate function is maximal and continuous are described. 
For rational two markings the index of the Veech group in $SL_2(\Z)$ is 
calculated in two different ways. 
\end{abstract}
\tableofcontents

\section{Introduction}
The motivation of this paper is the main Theorem in 
``Pointwise asymptotic formulas on flat surfaces" of A. Eskin and H. Masur
\cite{esmr98}. They prove that on strata of quadratic differentials $QD(g,P)$,
parametrizing half translation surfaces, there exists almost everywhere a limit 
of the quadratic growth rate and it is constant 
 (with respect to a natural $SL_2(\R)$ invariant Liouville measure $\mu$ on $QD(g,P)$
 see Veech \cite{ve86}).
The following questions arise:\\
1. Does the quadratic growth rate  associated to these constants majorize all other 
quadratic growth rates?\\
2. How can one compute the constants?\\
3a. Are there (sub) spaces of quadratic differentials, where one can prove not only
almost everywhere results but pointwise results?\\
3b. And if one has pointwise results, on which sets do the constants take ``different" 
values?\\ 

Our example the simplest one can find, but it has the advantage that one
can compute things elementary. The above stated questions could all be answered 
in a precise way.\\
The paper contains results about the growth rates of cylinders of periodic
trajectories or saddle connections on the two dimensional torus $\T^2:=(\R^2/\Z^2,dz)$, $dz=dx+idy$
with $n$ marked points. If the number of cylinders of periodic orbits (po)
or saddle connections (sc) of length less than $T$ on $\T^2_{[x_0,...,x_{n-1}]}$
is denoted by $N_{po/sc}(\T^2_{[x_0,...,x_{n-1}]},T)$ then the existence
of the limits
\[\lim_{T\rightarrow \infty}\frac{N_{po/sc}(\T^2_{[x_0,...,x_{n-1}]},T)}{T^2}\]
for every constellation of marked points i.e. for all $[x_0,...,x_{n-1}]\in {\T^{2n}}$ is proved.
\footnote{For all arguments the Parameterspace could be restricted to the case $x_0=0$, because this 
differs from the general marking only by a translation. All our counting results are invarinat under 
translations of the marking.}
 
Remark: In the mathematical sense this is not the true parameterspace, because we would like to have 
all marked points different, but since we make either pointwise statements with respect to one marking
or only measure theoretic statements with respect to the Lebesgue measure on this space
the difference to a well choosen parameter space plays no role. Moreover all statements are
true for any conformal class the tori might have. 
 
The existence of the limits is established by counting lattice
points of certain lattices in the universal cover $\R^2$ of $\T^2$. As a biproduct
of this method we conclude that there is a subset  $\T^{2n}_{gp(sc/po)} \subset \T^{2n}$ 
(which depends on what one counts: saddle connections or periodic families) 
of measure one in the parameter space of all n markings 
(contained in the set of all non rational n markings) 
where the above limit takes a maximal value $c^{max}_{sc/po}$. Moreover it is easy 
to see that if one approaches a point of the set $\T^{2n}_{gp(sc/po)}$ the asymptotic constants
of the sequence converge to $c^{max}_{sc/po}$. In other words the map
\[[x_0,...,x_{n-1}] \longmapsto \lim_{T\rightarrow \infty}
\frac{N_{po/sc}(\T^2_{[x_0,...,x_{n-1}]},T)}{T^2}\]
is continous at the points of the set $\T^{2n}_{gp(sc/po)} $. In the case of two marked points
and only in this case the set $\T^{2}_{gp(sc/po)} $ does not depend on what is counted (po's or sc's). 
Without restrictions one can assume that the first marking is $[0]$, then 
the set of points in general position consists exactly of all the non rational points 
$\R^2\backslash \Q^2 \bmod \Z^2$. 
In this case the quadratic growth rate function for saddle connections, assuming $\gcd(p_1,p_2,n)=1$, is given by:
\begin{eqnarray}
\left[\frac{p_1}{n}, \frac{p_2}{n} \right] &\longmapsto& \frac{6}{\pi}\left
(1 + \prod_{p| n\ prime}\left(1-\frac{1}{p^2}\right)^{-1}
\sum_{\gcd(i,n)=1}\left(\frac{1}{i^2}-\frac{1}{n^2}
\right)\right) \quad \mbox{on } \Q^2 \nonumber\\ 
\left[x,y \right] &\longmapsto& \frac{6}{\pi}+\pi \quad \mbox{for } (x,y) \notin \Q^2
\end{eqnarray} 
The continuity at non rational points implies that the above formula approaches
$\frac{6}{\pi}+\pi$ as $(\frac{p_1}{n}, \frac{p_2}{n})$ approximates 
a non rational number. In contrast to the general for two marked points 
this continuity is calculated directly by an estimate.
Furthermore for the two marked torus we write down the Veech and affine groups for rational markings.
An argument shows that every Veech surface coming from a rational marking is 
a $SL_2(\Z)$ deformation of one of the family $F=\{\T^2_{[\frac{1}{n},0]}: \ n \in \N\}$.
Therefore all the Veech groups of rationally marked tori are isomorphic 
to one of the Veech groups $V_{[\frac{1}{n},0]}$ associated to the family of surfaces $F$.
Finally the asymptotic formula of Veech and the asymptotic counting of lattice points 
in the plane is used to compute the index $[SL_2(\Z):V_{[\frac{1}{n},0]}]$.
With the help of the description of the Veech groups for rational 2 markings we 
rediscover the index of some classical congruence subgroups of $SL_2(\Z)$.

It is now well known \cite{ems} that one can also use Ratners Theorem to obtain pointwise results 
in modulispaces of torus coverings branched over n marked points. 
Of course from the point of pure existence of the growth rates
this is more general and in fact covers the examples described in the paper. 
On the other hand the results in this paper are explicit and the extremality, or more restrictive 
the maximality of the almost everywhere constants are proved for the above cases.
In this direction one of the main objects of the paper was to see if there is 
something in between purely rational markings and markings in general position which has not
an extremal growth rate. This question also could be answered with ``yes" (see Theorem 
\ref{tlimes} and the preparatory Lemma \ref{schnigi}).
Unfortunately our elementary methods can not directly taken over to prove 
something on families of branched coverings of tori, but nevertheless with the 
more sophisticated approach in \cite{ems} one can compute all the Siegel Veech constants 
for two marked torus coverings which are not Veech. The results in \cite{ems} show, 
that in this cases the Siegel Veech constants, as functions of the degree 
of the covering, behave monotonically, too. 
In the light of Eskins and this results it would be of interest to compute and prove 
related statements for marked Veech surfaces and their branched coverings.

Acknowledgements:\\
I have to thank for the  hospitality of the CPT in Marseille Luminy in person
of Ricardo Lima and his group during the jear '99 
as well as Alex Eskin and Pascal Hubert for fruitful discussions. 
I owe thanks to Eugene Gutkin for pointing out an unprecise statement to me 
and the referee for finding a mistake. 
This paper is part of my thesis and I like to give many thanks to 
Tyll Kr\"uger and Serge Troubetzkoy for their continuous encouragement 
making this work possible. 

\section{Notation and preparation}


For $[x_0],...,[x_n]\in \R^2/ \Z^2$ let $\T^2_{([x_0,...,x_{n-1}])} $ be the translation torus
$\T^2=(\R^2/\Z^2, dz=dx+idy)$ marked at the points $[x_0],[x_1],...,[x_{n-1}]\in \T^2$ .
If we do not write down the zeroth entry $[x_0]$ in the index vector we assume $[x_0]=0+\Z^2$.
Our calculations are restricted to the torus which is isometric to
$(\R^2/\Z^2,g:=d\!x^2+d\!y^2)$ with the induced Euclidean metric. 
With respect to this choice the periodic directions (of the geodesic flow)
are exactly the ones with rational slope.
The results are not restricted to these special translation tori, they hold for all marked translation tori.
In fact all the arguments are invariant under $SL_2(\R)$ deformations of the marked tori and any 
marked translation torus is an $SL_2(\R)$ deformation of the studied ones. 
The vector $v_s$ {\bf associated} to a saddle connection $s$ 
(or a prime periodic orbit) is the unique vector $v_s$ in $(\R^2, d\!z)$ that 
maps isometrically to $s$ in $\R^2 \rightarrow \T^2$ under the natural projection map. 
Recall that in the case of a periodic orbit $s$ the associated vector $v_s$ characterizes
a whole cylinder of closed geodesics.
A rational two vector like $(\frac{p}{n},\frac{q}{n})$ is always assumed to fulfill the condition 
$\gcd(p,q,n)=1$. 
\begin{defin} \nopagebreak $\mbox{}$ \linebreak
\begin{itemize}
\item[1.] Two points $[x],[y]\in \T^2$ are said to be {\bf relatively rational}, if
$x-y \in \Q^2 $.
 \item[2.]$SC_{(x,y)}$ denotes the set of saddle connections between the points
$[x]\in\T^2$ and $[y]\in\T^2$.
\item[3.] $l(s)$ is the length of a $s\in SC_{(x,y)}$ and finally if 
$[w],[z]\in s \in SC_{(x,y)}$ let  $l_{s}(w,z)$ be the distance between $[w]$  
 and $[z]$ \bf{measured along $s$}. 
\end{itemize}
\end{defin}
Relatively rational is obviously an equivalence relation. Every marking splits up in classes
of relatively rational points. Each class defines a marking, eventually of lower order.
By the results of Veech \cite{ve89} 
or Gutkin and Judge \cite{gj97}  $\T^2_{[x_0,...,x_{n-1}]}$ 
is a covering of the one marked torus $\T^2_{[x_0]}$ if and only if all marked points 
$x_i$ are pairwise relatively rational. 
Thus $\T^2_{[x_0,...,x_{n-1}]}$ is a  Veech surface if and only if 
for all markings are in one and the same marking class 
and from \cite{gj97} it follows:
$V\left(\T^2_{[x_0,...,x_{n-1}]}\right)\subset SL_2(\Q)$, where 
$V\left(\T^2_{[x_0,...,x_{n-1}]}\right)$ denotes the Veech group of
$\T^2_{[x_0,...,x_{n-1}]}$. 
Any affine diffeomorphism $\phi$ of the unmarked torus $(\T^2,dz)$ is 
given by an element $(A,[v])\in SL_2(\Z)\oplus \T^2$: $\phi([z])=[Az+v]$. 
Since every affine diffeomorphism of the marked torus defines one of the unmarked torus
by forgetting marked points one obtains: 
\begin{equation}\label{sl2}
V\left(\T^2_{[x_0,...,x_{n-1}]}\right)\subset SL_2(\Z). 
\end{equation}


\begin{defin}
A {\bf lattice} $G$ is a subset of $\R^2$ of the form
\[G=(p_1+q_1\Z,p_2+q_2\Z)\quad p_i,q_i \in \R^+\]
A {\bf point distribution} is a finite union of lattices.
The set of {\bf visible points} $G^V$ of of a discrete subset $G \subset \R^2$ consists 
of all points  $p\in G$ such there is no $0<|\lambda| <1$ and a point 
$q\in G$ with $q=\lambda p$.
\end{defin}
In other words: There is no other point of $G$ 
on the line between the origin and the points $\pm p$.
For the correct counting of periodic orbits and saddle connections
we need further:
\begin{defin}
If a point distribution $G$ is contained in another point distribution
$H$, then the visible points $H^{V_G}$ of $H$ with respect to $G$
are those $p\in H^V$ for which there is a $\lambda \neq 0$ such that
$\lambda p \in G$.
$G \subset H$ is called ${\mathbf H}$ {\bf complete} if $H^{V_G}=G^V$.
\\
$G$ is called {\bf closed}, if $mG \subset G$ for any integer $m \in \Z \backslash \{0\}$
such that there exist two points $x,y \in G$ with $x=my$.
\end{defin}

Let $[x],[y]\in \T^2$, again we put $x=0$. Then the set of saddle connections 
$SC_{(x,y)}$ is described by the visible points of the lattice of associated vectors 
$G_{(x,y)}=y+\Z^2$.\\
For a discrete subset $G \subset \R^2$ (in general point distributions) we define for $T>0$  
\[N(G^{(V)},T):=card(G^{(V)}\cap B(0,T)).\]
Here $B(0,T)$ is the (closed) disc of radius $T$ centered in the origin of $\R^2$.

Given a nonrational (marked) point $[y]$ on the marked torus $\T_{[0]}^2$. 
Take a saddle connection $s_{y}\in S_{(y)}:=S_{(0,y)}$ from the origin to $[y]$. 
It might be represented in $\R^2$ by $s_{y}=y +(m_y,n_y)$ with two fixed integers 
$m_y,n_y$ and an $y$ choosen from the set $[0,1)^2$. 
Assume we take a second point $[x]\in \T_{[0]}^2$ $x \in [0,1)^2$, different from $[y]$ 
so, that there is a saddle connection $s_{x}\in S_{(x)}$, 
which fulfills an equation 
\begin{equation}\label{relation}
s_{x}=x +(m_{x},n_{x})=r s_{y},
\ \mbox{with }\   r \in (0,1) \mbox{ and } n_x, m_x \in \Z.
\end{equation}
This means after mark the point $[x]$ the saddle connection $s_{x}$ destroys 
$s_{y}$ and generates two new saddle connections, $s_{x}$ and one in the set $S_{(x,y)}$. 
It is clear by equation (\ref{relation}), that
if the factor $r$ is a rational number, say $r=p/q$, all saddle 
connections represented by $y +(m_y,n_y) + q(k,l)$, $k,l \in \Z$ are divided by the marking 
$x$ to give saddle connections in $S_{(x,y)}$ and $S_{(x)}$. 
Now assume $[x]$ rationally divides two saddle connections from $S_{(y)}$ with the equations:
\begin{eqnarray}
x +(m_1,n_1)=\frac{p_1}{q_1}(y+(k_1,l_1))\\
x +(m_2,n_2)=\frac{p_2}{q_2}(y+(k_2,l_2))
\end{eqnarray}
with integers $k_i,l_i,m_i,n_i,p_i,q_i$, $i=1,2$ and $x,y$ as above. 
Subtract the two equations to obtain: 
\[(m_1,n_1)-(m_2,n_2)= \left( \frac{p_1}{q_1}-\frac{p_2}{q_2} \right)y +\frac{p_1}{q_1}(k_1,l_1)-\frac{p_2}{q_2}(k_2,l_2)\]
Since $y$ is non rational $ \frac{p_1}{q_1}=\frac{p_2}{q_2}$. 
Assuming this we have: 
\[(m_1,n_1)-(m_2,n_2)= \frac{p_1}{q_1}((k_1,l_1)-(k_2,l_2)).\]
This equation has solutions if $(k_1,l_1)-(k_2,l_2)$ is (componentwise) divisible by $q_1$. 
Thus the saddle connections in the set $S(x)$ dividing the ones 
from $S(y)$ rationally, all have a common stretch factor $p_1/q_1$. Hence, if there is 
more than one rationally dividing saddle connection, then they are all  
given by the (lattice) set of vectors 
\[\left\{\frac{p_1}{q_1}(y+(k,l)):\ (k,l)\equiv (k_1,l_1)\mod{q_1} \right\} \] 

If $[y]$, $p_1,q_1$ are fixed the different sets of rationally dividing saddle connections 
for the stretch factor $p_1/q_1$ are parameterized by $(k_1,l_1) \in \{0,...,q_1 -1\}^2$.    
It is now clear that $S_{(x,y)}$ contains a lattice. 

\begin{figure}

\epsfxsize=7truecm

\centerline{\epsfbox{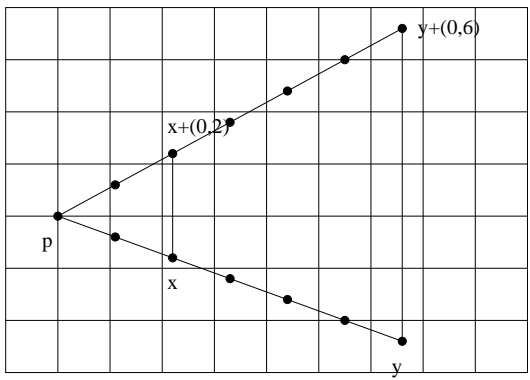}}
\caption{shows the `rational' intersection points of two saddle connections from $p$ to $y$ with the same horizontal
coordinates.}
\end{figure}
To make shure there are not too much 
parallel saddle connections in $S_{(x)}$ and $S_{(y)}$ having non rational length ratio, consider the set
$S^{nr}_{(x,y)} \subset S_{(x,y)}$ of saddle connections that are represented by $r\alpha$, where $r \in (0,1)$ is non rational. 

\begin{prop}\label{countnrat}
$N(S^{nr}_{(x,y)},T)= O(T)$. 
\end{prop}
\begin{proof}
Assume that the growth rate in question is of bigger than linear order, then there must be an 
integer $i \in \Z$ where one can find two saddle connections 
$s_{1},s_{2} \in S^{nr}_{(x,y)}$ represented by the points $s_1=r_1(y+(i,j_1))=r_1s_{1,y}$ and $s_2=r_2(y+(i,j_2))=r_2s_{2,y}$, 
$s_{1,y},s_{2,y}\in S_{(y)}$, $j_1,j_2 \in \Z$ (we assume $j_1<j_2$). 
If that would be not the case, for each $i \in \Z$ there could be at most one saddle connection of the 
form $r(y+ (i,l))$ in $S^{nr}_{(x,y)}$, therefore $N(S^{nr}_{(x,y)},T)\leq 2T$ contradicting our assumption. 
Now by the Strahlensatz of elementary geometry (see figure 1) each intersection point of the saddle connections  
$s_{1,y}$ and $s_{2,y}$ as above, must divide the length in rational amounts contradicting the assumption.  
\end{proof}

\begin{defin} Let $SC_{(x,y)}$ and $SC_{(x,z)}$ be two sets of saddle connections 
on the marked tori $\T^2_{[x,y]}$, $\T^2_{[x,z]}$ and suppose $[z]$ is 
rationally divided by $[y]$.
Then $G^{(x,z)}_{(x,y)}\subset G_{(x,y)}$ is the set of vectors associated to 
saddle connections $S_{(x,y)}$ on $\T^2_{[x,y,z]}$ which are parallel to vectors associated to 
saddle connections in $S_{(x,z)}$ and shorter in length.\\
\end{defin}

To summarize:
\begin{lemma}\label{schnigi}
Let $[x],[y],[z]\in \T^2$ be different points and $z-x \in [0,1)^2 \notin \Q^2$.
If $[y]\in s \in SC_{(x,z)}$, then $G^{(x,z)}_{(x,y)}$ contains a lattice, 
if and only if there is a rational number $p/q \in \Q$ such that 
$l_{s}(x,y)= \frac{p}{q} l_{s}(x,z)=\frac{p}{q} l(s)$, moreover the number $p/q$ is unique in the sense that 
this equation holds for any saddle connection represented by elements of $G^{(x,z)}_{(x,y)}$.
\end{lemma}
\begin{proof}
These are just the statements above translated to the situaion $[x]\neq [0]$. 
\end{proof}

The Lemma is viewed as the generalisation of the following observation on the torus 
$\T^2=\R^2/\Z^2$ (with the origin as a marked point): a point $p$ is rational if and only if
it is in the intersection of saddle connections. If this is the case the 
saddle connections through $p$ give rise to a lattice (by taking representatives in $\R^2$). 

The moral of the Lemma is: when counting saddle connections for general markings,
all markings which are lying on some saddle connection and dividing the length 
of it rationally will cause a lattice of intersection points. 
This is important as we will see later, this decreases the asymptotic constant of that marking. 
Because of this observation we also have to take care about the definition of the set of
markings in $\T^{2n}$, where we assume the growth rates to be extremal.
Thus for saddle connections we define the set of markings `in general position' (shortly: `gp') 
in $\T^{2n}$ inductively:
\begin{equation}
\begin{array}{l}
\T^{4}_{gp(sc)}=\{[x_0,x_1]\in \T^{4}: \  x_1-x_0 \notin \Q^2 \}\\
\mbox{ } \\
\T^{2(n+1)}_{gp(sc)}=\left\{ 
\begin{array}{l}
[x_0,...,x_n]\in \T^{2(n+1)}:[x_0,...,x_{n-1}]\in \T^{2n}_{gp(sc)},[x_n] \in \T^2 \ \mbox{and}  \\
\forall \ i \in \{ 0,...,n \} \  x_i \mbox{ does not divide any saddle connection}\\ 
\mbox{defined by the marking }[x_0,...,x_{i-1},x_{i+1},...,x_{n}]
\mbox{ rationally}
\end{array} 
\right\}.
\end{array}
\end{equation}

Thus the set of non extremal n markings is a set of full Lebesque measure
on $\T^{2n}$ (see Theorem \ref{tlimes}). 
For periodic orbits the set of markings in general position is defined in the following way:
\begin{equation}
\T^{2n}_{gp(po)}=\{[x_1,...,x_n]:\  x_j-x_i \notin \Q \  \forall \  i,j \in \{0,...,n\} \mbox{  and } j \neq i \}. 
\end{equation}
Observe that for all $n \in \N$ $\T^{2n}_{gp(sc)}\subseteq \T^{2n}_{gp(po)}$. The sets are equal 
in the cases $n=1$ (trivial) and $n=2$, for more marked points the inclusion is strict. 
It will turned out the markings in general position are exactely the sets of markings on which the growth rates are maximal.
Note that in the case of two marked points the set of markings in general position is the same as the 
set of non rational markings, always assuming one of the markings to be $[\vec{0}]$.
It should be noted here, that the pure existence of the growth rate functions 
for any choice of marking could also be seen using the Siegel Veech Formula together with 
Ratners Theorem for ergodic measures on homogeneous spaces (see \cite{ems}). 
The more elementary argumentation here and in the chapter on general markings shows, 
that even for complicated markings one has to count only simple lattice 
point distributions to obtain the quadratic growth constant.  

\section{Quadratic growth rates of lattices and point distributions in $\R^2$}

To compute quadratic growth rates we recall some preparatory Lemmas.
The first is quite well known and uses an argument going back to Gauss. 
\begin{lemma}\label{gauss}
Let $G=\{ (x,y)\in \R^2: (x,y)=(p_1 +q_1 \Z,p_2 +q_2\Z) \}$.
Then \[N(G,T)= card(G\cap B(T,0)) = \frac{\pi}{q_1q_2}T^2 + o(T^2).\]
\end{lemma}
\begin{proof} 
To begin we treat the case $p_1=p_2=0$. By drawing horizontal and vertical lines
through the lattice points, we got a decomposition of the plane in rectangles with
sidelengthes $(q_1,q_2)$.  Every lattice point of $G$ can be viewed 
as the upper right vertice of one rectangle. So the number of points of $G$ in the 
concentric circle $B(0,T)$ of radius $T$ is determined by the number of rectangles 
in that circle up to a certain correction term coming from the fact that a rectangle might 
be only partly inside the ball. But after increasing the radius of the ball 
by the length of the diagonal of the rectangles $\sqrt{q^2_1+q^2_2}$ the estimate 
\[N(G,T)\leq \frac{Vol\left(B\left(0,T+\sqrt{q^2_1+q^2_2}\right)\right)}{q_1q_2}\]
holds. And by decreasing the radius with the same amount we have:
\[\frac{Vol\left(B\left(0,T-\sqrt{q^2_1+q^2_2}\right)\right)}{q_1q_2})\leq N(G,T)\]
Putting together these two estimates leds to:
\[ \left|N(G,T)-\frac{\pi}{q_1q_2}T^2\right|\leq \frac{\pi}{q_1q_2}
\left(2T\sqrt{q^2_1+q^2_2} +q^2_1+q^2_2\right)\]
In the general case we have the estimate:
\begin{eqnarray}
\left|(G-(p_1,p_2))\cap B\left(0,T-\sqrt{p^2_1+p^2_2}\right)
\right| \leq \left|G \cap
B(0,T)\right| \leq \nonumber \\
\leq \left|\left(G-(p_1,p_2)\right)\cap B\left(0,T+\sqrt{p^2_1+p^2_2}\right)\right|
\end{eqnarray}
the statement follows easily.
\end{proof}
Before stating the next proposition let us make a general remark on our strategy of counting.
If our interest is for example counting prime closed geodesics, we associate a certain point 
distribution in $G \subset \Z^2$ to all geodesics in question. 
This point distribution $G$ does not reflect 
a priory the prime geodesics, but it will be choosen in a way that its quadratic asymptotics 
as in the last lemma can be calculated easily. The points in $G$ which are associated to 
prime geodesics will be the visible points of $G$ with respect to $\Z^2$.
The price we have to pay is that we cannot relate the distribution of visible points and
all points of (the completed and closed) $G$ just by the factor $\frac{6}{\pi^2}$.
We have to choose our multiplicities with respect to $G$:  
\begin{prop}\label{vielf}
Let $G\subset\Z^2$ be a  $\Z^2$ complete and closed  lattice. 
Let $K\subset \Z \backslash \{0\}$ be the set of numbers $l$ such that
$lG^{V}\subset G$. 
Then the relation 
\begin{equation}
N(G^{V},T)=N(G,T)\left(\sum_{l \in K}\frac{1}{l^2}\right)^{-1} + o(T^2)
\end{equation}
 holds between the numbers of points of $G$ and the number of visible points $G^V$ of $G$. 
 For the lattice $G=\Z^2$ this implies:
\begin{equation}
 N(\Z^2,T)=N({\Z^2}^{V},T)\frac{3}{\pi^2} + o(T^2).
\end{equation}
\end{prop}
\begin{proof}
First let us check that there is no $ x \in G$ that is not an integer multiple of some
$y \in G^V$. Assume there is, then there exists a $z \in G^V$ and a rational number
$\frac{p}{q}$ with $x=\frac{p}{q}z$. So $q$ divides $z$, a contradiction to the 
$\Z^2$ completeness of $G$. Since $G$ is closed and every element in $G$ is
some $K$-multiple of an element in $G^V$, thus:
\begin{eqnarray} \label{factor}
N(G,T)= \sum_{l\in K}N\left(G^{V},\frac{T}{l}\right) +o(T^2) = 
\left(\sum_{l\in K}\frac{1}{l^2}\right)N(G^{V},T)+o(T^2).
\end{eqnarray} 
The lattice $G=\Z^2$ is obviously $\Z^2$ complete, closed and $K=\Z\backslash \{0 \}$, thus 
\begin{eqnarray}
N(\Z^2,T)&=&\left(\sum^{\infty}_{\substack{l=-\infty \\ l\neq 0}}
\frac{1}{l^2}\right)N({\Z^2}^{V},T)+o(T^2) = \nonumber \\
&=&2\left(\sum^{\infty}_{l=1}\frac{1}{l^2}\right)N({\Z^2}^{V},T)+o(T^2)=
\frac{\pi^2}{3}N({\Z^2}^{V},T)+o(T^2). \nonumber
\end{eqnarray}
\end{proof}
Counting only a certain percentage $p$ of the length of each line 
in a point distribution $G$  is equivalent to counting the length of 
vectors in the  point distribution \\ $pG:=\{pv: v \in G \}$. Thus there is a relation:
\begin{equation}\label{scaling}
N(pG,T) = N\left(G,\frac{T}{p} \right)= \frac{1}{p^2}N(G,T)+o(T^2).
\end{equation}
Remark: As is seen above the relating factor between cylinders of geodesics 
and prime geode\-sics will be $\frac{3}{\pi^2}$ instead of $\frac{6}{\pi^2}$ because we don't 
distinguish between the two directions a prime geodesic might have. Our decision to do that
is based on the point of view : count only (the length of) the geometric object.
\begin{defin}  Define the set of lattices parametrized by $ p_1,q_1,p_2,q_2 \in \Z $, with 
$\gcd(p_i,q_i)=1$ as:
\[G_{(q_1,p_1,q_2,p_2)}:=\left(q_1 \Z+ p_1,q_2 \Z+ p_2\right) \subset \Z^2.\]
Further for $n \in \N$ 
and a divisor $m$ of $n$ define $I^n_m$ to be the set of numbers which are 
relatively prime to $\frac{n}{m}$
\[I^n_m= \left\{j\in \left\{1,2,...,\frac{n}{m}-1\right\}: \gcd \left(j,\frac{n}{m}\right)=1\right\}\]
(if $m=1$ we write shortly $I^n$ instead of $I^n_m$) and the point distributions
\[G_{I^n_m}:=\bigcup_{k\in I^n_m}G_{(n,mk,n,0)}\subset \Z^2.\]
\end{defin}
With this conventions 
\begin{lemma}\label{pvert}
The point distributions $G_{I^n}$ are for all $n \in \N$ complete with respect to $\Z^2$
and closed. We have:
\begin{equation}\label{vf}
N(G_{I^n},T)=2N(G^{V}_{I^n},T)\left[\sum_{k\in I^n} 
\sum^{\infty}_{l=1 }\frac{1}{(nl+k)^2}\right].
\end{equation}
As a consequence, for all prime numbers $p$
\begin{equation} \label{primvf}
N(G_{I^p},T)=N(G^{V}_{I^p},T)\frac{\pi^2}{3}\left(1-\frac{1}{p^2}\right).
\end{equation} 
\end{lemma}
\begin{proof}
First one has to find all integers which map $G_{I^n}$ under elementwise
multiplication to  $G_{I^n}$ itself.
For $l \in \Z$ and $p\in G_{I^n}$ we have $lp\in G_{I^n_m}$, if and only if 
\[l \equiv r \mod{n} \quad , \mbox{with } r \in I^n.\]
Therefore $l$ has to be in $ n\Z+r \quad \mbox{with }\ r \in I^n$.
Obviously the multiplication with a number from this set maps points of $G_{I^n}$ again to
$G_{I^n}$. 
So if the point distributions are complete with respect to $\Z^2$ they are also 
closed. To see the last statement take $l\in \Z \backslash \{0\}$
and choose $(p_1,p_2)\in \Z^2$ with $l(p_1,p_2) \in G_{I^n}$.
This means $l(p_1,p_2) = (kn+r,mn)$ with $r \in I^n$ and $k,m \in \Z$.
The equation in the first coordinate says that $l$ and $p_1$ are relatively prime 
to $n$ and from this it follows that $p_2$ has to be a multiple of $n$. Thus $(p_1,p_2)$
is already contained in $G_{I^n}$ and the completeness follows.
Using equation (\ref{factor}) from Proposition \ref{vielf} 
gives the first part of the Lemma.
For the second part we simply observe that because $p$ is prime $I^p=\{1,2,...,p-1\}$:
\[\sum_{k\in I^p} \sum^{\infty}_{l=1 }\frac{1}{(pl+k)^2}=
\sum^{\infty}_{l=1 }\frac{1}{l^2}-\sum^{\infty}_{l=1 }\frac{1}{(pl)^2}\] 
\end{proof}
Remark: Obviously the lattices $I^{n}_m$ are for $m\neq 1$ never $\Z^2$ complete.
Before we compute the relating factor of equation (\ref{vf}) in general, an
observation is helpful: Let ${\chi }_{0}: \Z \rightarrow \Z/2\Z$ be the principal 
character modulo $n$, 
that is ${\chi}_0(l)=1$ if $\gcd(l,n)=1$ and is equal to $0$ if not, then obviously:
\[\sum_{k\in I^n} \sum^{\infty}_{l=1 }\frac{1}{(nl+k)^2}=
\sum^{\infty}_{l=1 }\frac{{\chi}_{0}(l)}{l^2}= L({\chi}_0,2)\] 
where $L({\chi}_0,2)$ is the $L$ series with respect to ${\chi}_0$ evaluated at $2$. 
The following is well known 
(see for example exercise (3.d) page 166 in \cite{ten} or Serre \cite{se}) 
\begin{prop} \label{sumform}
If $s$ is a complex number with real part bigger than $1$
then:
\begin{equation} 
L({\chi}_0,s)=\sum^{\infty}_{l=1 }\frac{\chi_{0}(l)}{l^s}=\zeta(s)\prod_{p|n\ prime}
\left(1-\frac{1}{p^s}\right)
\end{equation} 
\end{prop}
\begin{proof}
The condition on the real part of $s$ guaranties absolute convergence
of the series so that the following manipulations are justified.
We start with the right hand side:
\begin{eqnarray}
\zeta(s)\prod_{p|n\ prime}\left(1-\frac{1}{p^s}\right)=
\sum^{\infty}_{l=1}\frac{1}{l^s}\left(\sum_{m|n}\frac{\mu(m)}{m^s}\right)= \nonumber \\
=\sum^{\infty}_{k=1}\sum_{m|\gcd(k,n)}\frac{\mu(m)}{k^s} =
\sum^{\infty}_{k=1}\frac{1}{k^s}\left(\sum_{m|\gcd(k,n)}\mu(m)\right)=
\sum^{\infty}_{k=1 }\frac{{\chi}_{0}(k)}{k^s}
\end{eqnarray}
We have used the well known M\"obius $\mu$ function, the formula 
\[\sum_{d|n}\mu(d)=\left\{ \begin{array}{l}
1 \mbox{ if } n=1\\
0 \textrm{ if } n>1 
\end{array} \right.  \] 
and summation over $k=lm$ to reorder the sum in a useful way. 
\end{proof}
Evaluating the above $L$ series for $s=2$ we find
\begin{equation}
N(G_{I^n},T)=N(G^{V}_{I^n},T)\frac{\pi^2}{3}\prod_{p|n \  prime}
\left(1-\frac{1}{p^2}\right)+o(T^2).
\end{equation}
The left hand side of this equation is easy to compute with the help of 
Lemma \ref{gauss} and the definition of $G_{I^n}$ 
\begin{eqnarray}
N(G_{I^n},T)=\frac{\pi}{n^2}|I^n|T^2+o(T^2).
\end{eqnarray}
Use Eulers $\varphi$ function to identify
\begin{equation}\label{in}
|I^n|=\varphi(n)=n\prod_{p|n \ prime}\left(1-\frac{1}{p}\right)
\end{equation}
and finally write
\begin{eqnarray}\label{vffinal}
\lim_{T \rightarrow \infty}
\frac{N(G^{V}_{I^n},T)}{T^2}&=&\frac{3}{\pi}\frac{\varphi(n)}{n^2}\prod_{p|n \ prime}
\left(1-\frac{1}{p^2}\right)^{-1}=
\nonumber\\
&=&\frac{3}{n\pi}\prod_{p|n \ prime}\left(1+\frac{1}{p}\right)^{-1}.
\end{eqnarray}
For later use we show
\begin{prop}\label{completion}\hfill \\
$G_{I^n}$ containes the $\Z^2$ completion of the lattice 
$G_{(n,r,n,0)}=$$\{(n\Z+r,n\Z):\  r \in I^n \} $. 
\end{prop} 
\begin{proof}
It is enough to show that a multiple of each point $(kn+l,mn) \in G_{I^n}$
is in $G_{(n,r,n,0)}$. That is we have to solve $rx \equiv l \pmod{n}$.
But this is always possible since $r$ and $l$ are relatively prime with respect
to $n$. 
\end{proof}

\section{The Veech group of $\mathbf{\T^2_{\left[\frac{p_1}{q_1},\frac{p_2}{q_2}\right]}}$ and its  
index in $\mathbf{SL_2(\Z)}$}

We now describe the affine group $A\mathit{ff}(\T^2_{x})$ for a rational point $x\in \Q^2$. 
\begin{prop} \label{afftor}
Assume $(\frac{1}{2}, \frac{1}{2}),(0, \frac{1}{2}), (\frac{1}{2}, 0)\neq x\in \Q^2$.
Then 
\[A\mathit{ff}(\T^2_{x})=\left\{z \mapsto Az+v:\! A\in SL_2(\Z) \mbox{ with \,}\left\{\! \begin{array}{ll}
 Ax \equiv -x  \pmod{\Z^2} & \textrm{ if } v=x\\
  Ax  \equiv x   \pmod{\Z^2} & \textrm{ if } v=0 
\end{array} \right. \!\right\} \]
If $y=(\frac{1}{2}, \frac{1}{2}),(0, \frac{1}{2}), (\frac{1}{2}, 0)$ we have
\[ A\mathit{ff}(\T^2_{y})=\left\{z \mapsto Az+v: A\in SL_2(\Z) \mbox{ with }
\left\{\begin{array}{l}
 Ay \equiv y  \pmod{\Z^2} \\
 v=y 
\end{array} \right.\right\} \]
\end{prop}
\begin{proof}
$A\in SL_2(\Z)$ is a consequence of the discussion around equation (\ref{sl2}) in the preparation section. 
The only thing which is left to prove is the statement on the translation vector $v$. 
It is defined only modulo $\Z^2$, thus one can assume $v\in[0,1)^2$.
The subgroup $N(\T^2_x)$ of pure translations is trivial if 
$x\neq (\frac{1}{2}, \frac{1}{2}),(0, \frac{1}{2}), (\frac{1}{2}, 0)$:\\ 
Assume the lattice $\Z^2_x=x+\Z^2$ generated is by use of $v$ translated to $\Z^2$ and 
by the same translation $\Z^2$ is moved to $\Z^2_x$, then $v\equiv x \equiv -x   \pmod{\Z^2}$ must hold. 
Thus either $v$ is trivial, or $x$ is one of the exceptions. 
In general the same problem comes up if the map should interchange the lattices 
$\Z^2$ and $\Z^2_x$. The only possibility is: $A$ has to map $\Z^2_x$ to
$\Z^2_{-x}$ and $v$ translates this into $\Z^2$. Thus we obtain the conditions:
\[Ax \equiv -x   \pmod{\Z^2}\ \mbox{ and } v=x \quad \mbox{or}\quad  Ax\equiv x
 \pmod{\Z^2}\ \mbox{ and }v=0\]
For the three exceptions there is in fact only one condition,
because in these cases $-x \equiv x  \pmod{\Z^2}$, thus $-id$ is automatically contained 
in these affine groups. 
Moreover the translation $v=x$ is contained in the groups.
\end{proof}
\begin{satz} \label{kong}The Veech group $V(\T^2_x)$ where  
$x=\left[\frac{p_1}{q_1},\frac{p_2}{q_2}\right]$ 
\begin{itemize}
\item for $ 1\leq \gcd(q_1,q_2)< \min(q_1,q_2) $ is
\begin{displaymath}
V(\T^2_x)=
\left\{\begin{array}{cc} \pm \left( 
\begin{array}{cc} a & b \\ c & d                        
\end{array} \right) \in SL_2(\Z):&
                    \begin{array}{cc}
                         a\equiv 1 \quad \pmod{q_1} & b \equiv 0  \pmod{q_2} \\
                         c\equiv 0 \quad \pmod{q_1} & d  \equiv 1 \pmod{q_2}
                        \end{array} 
        \end{array}
\right\}
\end{displaymath}
\item for $ q_1|q_2 $ (or  $ q_2|q_1 $ analogously) is
\begin{displaymath}
V(\T^2_x)=
\left\{ \begin{array}{cc} \pm \left( 
                     \begin{array}{cc}
                         a & b \\
                         c & d 
                        \end{array} \right) \in SL_2(\Z): &
                         \begin{array}{c}
                         a+\frac{q_1p_2}{q_2}b \equiv 1 \pmod{q_1} \\
                         b \equiv 0 \pmod{\frac{q_2}{q_1}} \\
                         c\equiv 0 \pmod{q_1} \\
                         d  \equiv 1 \pmod{q_2}
                        \end{array} 
        \end{array}
\right\}
\end{displaymath}
\item for $ q_1=q_2 $ is
\begin{displaymath}
V(\T^2_x)=
\left\{\!\begin{array}{ll} \pm  \left( 
                     \begin{array}{cc}
                         a & b \\
                         c & d 
                        \end{array} \right) \!\in SL_2(\Z):\!\! &
                         \begin{array}{c}
                          \left( 
                     \begin{array}{cc}
                         a & b \\
                         c & d 
                        \end{array} \right)
                        \left(\begin{array}{c} p_1 \\ p_2 \end{array}\right)
                                                              \equiv \left(\begin{array}{cc}
                         1& 0 \\
                         0 & 1 
                        \end{array} \right)
                                              \pmod{q_1}
                        \end{array} 
        \end{array}
\! \right\}
\end{displaymath}
In the case $q_1=2$ that is $(\frac{p_1}{q_1},\frac{p_2}{q_2})=(\frac{1}{2},\frac{1}{2})$ we 
have to add the rotation:
\begin{equation}\nonumber 
r_{\frac{\pi}{2}}=\begin{pmatrix} 
 _{0} & _{1} \\ 
 _{-1} & _{0}
\end{pmatrix}
\end{equation}
\item for $ q_2=0 $ (or equivalently  $q_1=0$) is
\begin{displaymath}
V(\T^2_x)=
\left\{\begin{array}{ll} \pm \left( 
                     \begin{array}{cc}
                         a & b \\
                         c & d 
                        \end{array} \right) \in SL_2(\Z):&
                         \begin{array}{cc}
                        \left( 
                     \begin{array}{cc}
                         a & b \\
                         c & d 
                        \end{array} \right)& \equiv \left( 
                     \begin{array}{cc}
                         1& b \\
                         0 & 1 
                        \end{array} \right) \pmod{q_1}
                        \end{array} 
        \end{array}
\right\}
\end{displaymath}
\end{itemize}
\end{satz}
\begin{proof}
By the last proposition we have to evaluate the condition 
\[ \begin{array}{cc}
                        \left( 
                     \begin{array}{cc}
                         a & b \\
                         c & d 
                        \end{array} \right) \left( 
                     \begin{array}{c}
                         \frac{p_1}{q_1}\\
                          \frac{p_2}{q_2}
                        \end{array} \right) \equiv
                      \pm\left(\begin{array}{c}
                         \frac{p_1}{q_1}\\
                          \frac{p_2}{q_2}
                        \end{array} \right) \pmod{\Z^2} 
        \end{array}
\] on the coefficents  $a,b,c,d \in \Z$ . If for example $p_2=0$
then we have
\[ \begin{array}{l}
                        \left( 
                     \begin{array}{c}
                         \frac{p_1}{q_1}a   \\
                         \frac{p_1}{q_1}c  
                        \end{array} \right) \equiv \pm \left( 
                     \begin{array}{c}
                         \frac{p_1}{q_1}\\
                          0
                        \end{array} \right) \pmod{\Z^2}
        \end{array}
\]
It follows directly that $a\equiv \pm 1  \pmod{q_1}$, $c \equiv 0  \pmod{q_1}$ and
$b\in \Z$ is arbitrary. Since the determinant of the matrix is always $1$ we have the condition 
$ad \equiv 1  \pmod{q_1}$. Therefore $d \equiv \pm 1  \pmod{q_1}$, if $a \equiv \pm 1  \pmod{q_1}$.
The other three cases are equally simple, we omit them.
\end{proof}
We prove now that all 2 marked tori are affine equivalent to one from a 
special family, if the marking is rational. 
\begin{prop}\label{reduction} If $x=\left[\frac{p}{n},\frac{q}{n}\right]$ ($\gcd(p_1,p_2,n)=1$), 
$\T^2_x$ is affine isomorphic to $\T^2_{\left[\frac{1}{n},0\right]}$. 
More precisely: there exists an affine diffeomorphism of $\T^2$  represented by an element of $SL_2(\Z)$ 
which maps $\left[\frac{p}{n},\frac{q}{n}\right]$ to $\left[\frac{1}{n},0\right]$
(and preserves $[0]$). 
\end{prop} 
\begin{proof}
The Veech group $V(\T^2)$ of the one marked torus is represented by $SL_2(\Z)$. 
Thus it is enough to show: given a point $[p/n,q/n]\in\T^2$ 
with $\gcd(p,q,n)=1$, there is a linear map in $SL_2(\Z)$ that maps it to the point 
$[1/n,0]\in\T^2$. 
Let $c:=\gcd(p,q)$ then the equation $\frac{q}{c}p-\frac{p}{c}q=0$ holds. We want 
the vector  $(-\frac{q}{c},\frac{p}{c})$ to become the second row of a matrix $A\in SL_2(\R)$, thus 
the vector $(a,b)$ representing the first row has to fulfill $1=\det(A)=a\frac{p}{c}+b\frac{q}{c}$.
This equation has  solutions  $a,b \in \Z$ because  $\frac{p}{c}$ and $\frac{q}{c}$ are relatively prime.
Taking solutions $a,b$ defines an $A\in SL_2(\R)$ with $A(p,q)=(c,0)$ with $c$ relatively prime to $n$.
Hence the congruence $kc\equiv 1 \pmod{n}$ has a solution $k \in \Z$ and one can use the linear Dehn twists 
\[\begin{array}{cc}
                        \left( 
                     \begin{array}{cc}
                         1 & 0 \\
                         k & 1 
                        \end{array} \right) \left( 
                     \begin{array}{c}
                         c\\
                          0
                        \end{array} \right) \equiv
                      \left(\begin{array}{c}
                         c\\
                          1
                        \end{array} \right)  
        \end{array} \pmod{n}\]
  and      
 \[\begin{array}{cc}
                        \left( 
                     \begin{array}{cc}
                         1 & -c \\
                         0 & 1 
                        \end{array} \right) \left( 
                     \begin{array}{c}
                         c\\
                         1
                        \end{array} \right) \equiv
                      \left(\begin{array}{c}
                         0\\
                          1
                        \end{array} \right)  
        \end{array}\pmod{n} \]
 to obtain the result up to a rotation of 90 degrees.       
\end{proof}
\begin{kor}\label{isom} 
The Veech groups $V(\T^2_x) \subset SL_2(\Z)$ for points 
$x=[\frac{p}{n},\frac{q}{n}]$ with $\gcd(p,q,n)=1$
are all isomorphic. The isomorphism is given by conjugation with an element of
$SL_2(\Z)$. Moreover the asymptiotic growth rates are the same, if the surfaces are 
affine isomorphic.
\end{kor}
\begin{proof}
The first two statements follow from the above proposition. Veechs asymptotic formula depends 
only on the conjugacy class of the Veech group in $SL_2(\R)$ and on the areas of the 
cylinders of closed geodesics (see for example Vorobetz paper \cite{va}), 
which are invariant under affine isomorphisms. In other words the asymptotic constants 
depend only on the $SL_2(\R)$ orbit of the given Veech surface in the moduli space. 
But Proposition (\ref{reduction}) shows the marked surfaces in question are all on one $SL_2(\R)$ orbit. 
\end{proof}
\begin{rem} Since every rational number $x$ has a unique representation as above, 
one can restrict the considerations to families of markings defined by 
\[x_{n}:=\left[\frac{1}{n},0\right]\quad   \mbox{or}
\quad  x_{(n,n)}:=\left[\frac{1}{n},\frac{1}{n}\right]\quad \mbox{where }\   1\neq n \in \N \]
and to the corresponding tori. We choose the family $\T^2_{x_n}$. 
\end{rem}
To compute the index $[SL_2(\Z):V(\T^2_{[\frac{p}{n},\frac{q}{n}]})]$, the first observation is that 
the only affine map $\phi \in A\mathit{ff}(\T^2_{[\frac{p}{n},\frac{q}{n}]})$ ($\gcd(p,q,n)=1$) 
which interchanges $[0]$ and $[\frac{p}{n},\frac{q}{n}]$ is given by rotation of 180 degrees 
and a translation by the vector $(\frac{p}{n},\frac{q}{n})$ (for $n>2$). Since $\phi$ maps 
$[\frac{p}{n},\frac{q}{n}]$ to $0$ it is of no interest for the following discussion. 
Up to the subgroup of generated by the idempotent element $\phi$, 
$A\mathit{ff}(\T^2_{[\frac{p}{n},\frac{q}{n}]})$ can be viewed as the isotropy subgroup  
$A\mathit{ff}_{[\frac{p}{n},\frac{q}{n}]}(\T^2)\subset A\mathit{ff}(\T^2) $. For a translation structure 
$u:=(S_{g},\omega)$ and $x\in S_{g}$ one defines the isotropy group $A\mathit{ff}_x(u)$ of $x$ 
and its associated Veech group in the usual way: 
\begin{eqnarray} &A\mathit{ff}_x(u) := \{g \in A\mathit{ff}(u): gx=x\} \quad \mbox{and}\nonumber \\ 
&V_x(u):=d\circ A\mathit{ff}_x(u)\subset V(u).\nonumber
\end{eqnarray} 
These are subgroups of $A\mathit{ff}(u)$, $V(u)$ respectively. 
As a direct consequence from the definitions and elementary group theory one obtains: 
\begin{prop}\label{shindex}
Let $u$ be a translation structure and $u_x$ the structure $u$ where the point $x$ is marked. 
Then 
\begin{equation}
[A\mathit{ff}(u):A\mathit{ff}_x(u)]= \ord(A\mathit{ff}(u)/ A\mathit{ff}_x(u))=
|O_{A\mathit{ff}(u)}(x)|
\end{equation}
and
\begin{equation}\label{vindex}
[V(u):V_x(u)]=\ord (V(u)/V_x(u))
\end{equation}
where $O_{A\mathit{ff}(u)}(x):=\{y\in S_{g}: y \in A\mathit{ff}(u).x\}$ is the orbit of $x$ 
under $A\mathit{ff}(u)$ and $G/H$ denotes the set of right cosets $xH$ of the subgroup $H\subset G$.  
\end{prop}
\begin{beweis}
See for example Chapter I \S 5 in S. Langs book on algebra \cite{langa}.
\end{beweis} 
Lemma \ref{reduction} above implies together with Proposition 
\ref{shindex} for $n>2$
\begin{eqnarray}[A\mathit{ff}_{[\frac{p}{n},\frac{q}{n}]}(\T^2):A\mathit{ff}(\T^2)]
&=&|\{(a,b): 0< a,b \leq n \mbox{ and } \gcd (a,b,n)=1 \}|\nonumber \\
&=& n^2\prod_{p|n \  prime}\left(1-\frac{1}{p^2}\right) 
\end{eqnarray}
The last equation on the cardinality of the set of pairs $(a,b)$ is well known. 
$d\phi=-id$ is normal in $SL_2(\Z)$ and therefore also in $V(\T^2_{[\frac{p}{n},\frac{q}{n}]})$, thus
\[V_{[\frac{p}{n},\frac{q}{n}]}(\T^2)\cong V(\T^2_{[\frac{p}{n},\frac{q}{n}]})/(-id).\] 
Since $-id$ generates a subgroup of order $2$ and the subgroup $N(\T^2_{[\frac{p}{n},\frac{q}{n}]}) \subset 
A\mathit{ff}(\T^2_{[\frac{p}{n},\frac{q}{n}]})$ of pure translations is trivial if $n\ne 2$ 
we have 
\[[SL_2(\Z):V(\T^2_{[\frac{p}{n},\frac{q}{n}]})]=\frac{n^2}{2}\prod_{p|n \  prime}
\left(1-\frac{1}{p^2}\right)\]
For $n=2$ the index is $3$ because $-id$ is already an element in 
$V_{[\frac{p}{n},\frac{q}{n}]}(\T^2)$.

In contrast to this calculation there is a second way to compute the index of the Veech groups 
using the asymptotic formula of Veech and the results from the last section. 
To begin, the Veech groups are presented as subgroups of $SL_2(\Z)$, but since we do not want 
to count orbits in both directions we take the quotient modulo $\pm id$: $PSL_2(\Z)$. 
Since both modulo $-id$ equivalent elements are in $V(\T^2_{x_n})$ we can do the same with 
the Veech groups and loose no information when taking the preimage in $SL_2(\R)$. 
Thus, denote by $V_{p}(\T^2_{x_n})$ the image of $V(\T^2_{x_n})$ in  $PSL_2(\R)$.
The orbits $V_{p}(\T^2_{x_{n}})(\pm v_{\CC_i})\in \R^2 ,\ i=1,...,k(x_{n})$ 
($k(x_{n})$ is the number of different orbits) of the vectors $\pm v_{\CC_i}$ 
associated to a maximal cylinder $\CC_i$ of closed geodesics which fills the torus under 
 $V_{p}(\T^2_{x_{n}})$ give always the same asymptotic constant 
 (for the formula see Proposition 6.3 in \cite{gj97}:
\[\lim_{T \rightarrow \infty}\frac{N(V_{p}(\T^2_{x_{n}})(\pm v_{\CC_i}),T)}
{T^2}=vol\left(\H/V_{p}(\T^2_{x_{n}})\right)^{-1}
\frac{[A\mathit{ff}(\CC_i):A\mathit{ff}_0(\CC_i)]}{area(\CC_i)}
\] 
Since each cylinder $\CC_i$ fills the torus we have $area(\CC_i)=1$ and 
$[A\mathit{ff}(\CC_i):A\mathit{ff}_0(\CC_i)]=1$.
Here $A\mathit{ff}_0(\CC_i)$ is the group generated by linear Dehn twists around $\CC_i$. 
Thus we can write:
\begin{eqnarray}
vol\left(\H/V_{p}(\T^2_{x_{n}})\right)&=&[PSL_2(\Z):V_{p}(\T^2_{x_{n}})] \  
vol\left(\H/PSL_2(\Z)\right)= \nonumber \\
&=&\frac{\pi}{3}[SL_2(\Z):V(\T^2_{x_{n}})].
\end{eqnarray}
The other side of Veech's asymptotic formula is the quadratic growth rate of the lengths of 
the maximal cylinders that completely fill the torus. Their directions are exactly the directions of lines
in $\R^2$ which begin at the origin and cross some point of the lattice $x_{n}+\Z^2$.
To get a set of integer points which represent (eventually multiples of the length of) 
the closed cylinders in question we simply multiply with $n$ and get the 
lattice $G_{(n,1,n,0)}=n\Z^2+(1,0) \subset \Z^2 $ studied earlier.  Since we are on the torus $\R^2/\Z^2$ 
to count the length spectrum of the cylinders correctly one has to look for the 
visible points of the completion of $G_{(n,1,n,0)}$ with respect to $\Z^2$. 
This completion, by Proposition \ref{completion}, is the point distribution $G_{I^n}$, 
thus finally: 
\begin{equation}\label{kxn}
\lim_{T\rightarrow \infty}\frac{N(G^{V}_{I^n}, T)}{T^2}=
\frac{3k(x_{n})}{\pi}[SL_2(\Z):V(\T^2_{x_{n}})]^{-1}
\end{equation}
What is left to calculate is $k(x_{n})$ : 
\begin{prop}\label{invariante}
\begin{equation}
k(x_{n})=\left\{ \begin{array}{lll}
\frac{1}{2}\varphi(n)=\frac{n}{2}\prod_{p|n \  prime}\left(1-\frac{1}{p}\right)  
&\mbox{ if}&  n>2 \\
 1 & \mbox{ if} &  n=2
\end{array}\right.
\end{equation}
\end{prop}
\begin{proof}
Observe that the ratio of the length of the two saddle connections connecting
$[0]$ and $[x_n]$ and bounding the same cylinder of periodic orbits is an invariant
under the operation of the Veech group. Furthermore modulo $V(\T^2_{x_{n}})$ we
can represent every periodic orbit or saddle connection by a point in 
$Q_n=\{(k,l) \in \Z^2 :\ 0< k,l\leq n\} $. Thus the directions with only one 
periodic family are represented by:
\[\{(i,n)\in Q_n: \gcd(i,n)=1\}\]
The line from $(0,0)$ to $(i,n)$ intersects exactly one point of the set $x_{n}+\Z^2$ 
characterized by  
\[\left(x\frac{i}{n},x \right)\  \mbox{ with } xi\equiv 1 \pmod{n}.\]
If $i\in \{1,2,...,n-1\}$ runs through the numbers relatively prime to $n$, then  
$x$ does this as well. Thus the above regarded length ratio is in any case one of the numbers
\[\frac{\min(x,n-x)}{\max(x,n-x)} \quad \mbox{ with }x\in I^n .\]
We  see, exactly two relatively prime numbers modulo $n$, namely $x$
and $n-x$ having the same invariant. If this invariant is complete, which will 
be shown in the next Lemma, then the number of different one cylinder orbits under the 
operation of the Veech group is given by
\[\frac{1}{2}\varphi(n)=\frac{1}{2}\left|\{x\in\{1,2,...,n-1\}:
\ \gcd(x,n)=1\}\right|=\frac{n}{2}\prod_{p|n\ prime}\left(1-\frac{1}{p}\right).\]
The exceptional case $n=2$ is trivial, because there exists only one direction of one cylinder
in $Q_2$. 
\end{proof}
\begin{lemma} Let $\T^2_x$ be a rational marked torus. Two directions $v_1,v_2$ with only one 
periodic cylinder are  in the same orbit under the operation of the Veech group, if they
have the same ratio  
\[I(v_i)=\frac{\min(s^{v_i}_1,s^{v_i}_2)}{\max(s^{v_i}_1,s^{v_i}_2)} \quad i=1,2\] 
Here $s^{v_i}_j\ j=1,2$ are the two saddle connections which are on the boundary of the 
cylinder in direction $v_i$.
\end{lemma}
\begin{proof}
Obviously the condition is necessary. All that is left is to show that it is sufficent.
That is we have to find an element of the Veech group that is a map between 
the two directions with the same invariant. Forgetting for a moment the marked
point $x$  we can map any periodic direction to any other by elements of $SL_2(\Z)$.
Especially the saddle connections in question are mapped on one another. Because 
the invariant is the same for both possible orientations of the saddle connections 
and  $-id\in SL_2(\Z)$ we have found a map, which also transports one marked point 
to the other.
\end{proof}
Now we can compute the index with the help of the counting formula again:
\begin{kor} Let $x=[\frac{p'}{n},\frac{q}{n}]$, $n>2$, $gcd(p',q,n)=1$ then
\begin{equation}
 [ SL_2(\Z):V(\T^2_{x})]=\frac{n^2}{2}\prod_{p|n \  prime}\left(1-\frac{1}{p^2}\right)
\end{equation}
holds. In the case $n=2$ the index is $3$.
\end{kor}
\begin{proof}
From equality (\ref{kxn}) together with $k(x_n)=\frac{1}{2}\varphi(n)$
for $n>2$ it follows 
\begin{equation}\label{index}
[SL_2(\Z):V(\T^2_x)]=\frac{3}{2\pi} \varphi(n)\left[\lim_{T\rightarrow \infty}
\frac{N(G^{V}_{I^n}, T)}{T^2}
\right]^{-1}.
\end{equation}
Then put in the right hand side of equation (\ref{vffinal}) to get the result.
If $n=2$ then $k(x_2)=1$ and we obtain $3$.
\end{proof}
\begin{rem}
The two ways of computing the index of the Veech groups of two marked tori might be used 
to obtain the asymptotic formula (\ref{vffinal}), counting the lattice points, without using the 
second section. One can take the first index computation and the Veech formula together 
with the knowledge about the number of cusps $k(x_n)$ (see Proposition \ref{invariante}) 
to compute $ \lim_{T\rightarrow \infty}N(G^{V}_{I^n}, T)/T^2$ backwards.
\end{rem}
Nevertheless the first given way to compute the index of the Veech groups seems to be a 
simple method to calculate for the index of the well known groups 
(for example compare S.Lang \cite{lang}
\begin{displaymath}
\Gamma_1(n)=\left\{\begin{array}{ll} \left( 
                     \begin{array}{cc} 
                         a & b \\
                         c & d 
                        \end{array} \right) \in SL_2(\Z):&
                         \begin{array}{cc}
                        \left( 
                     \begin{array}{cc}
                         a & b \\
                         c & d 
                        \end{array} \right)& \equiv \left( 
                     \begin{array}{cc}
                         1& b \\
                         0 & 1 
                        \end{array} \right)\pmod{n}
                        \end{array} 
        \end{array}
\right\}
\end{displaymath}
in $SL_2(\Z)$. Either the $\Gamma_1(n)$ are subgroups of index two in $V(\T^2_{x_n})$ if $n>2$, or 
in the case $n=2$ $-id \in \Gamma_1(2)$ and $V(\T^2_{x_2})$ is isomorphic to $\Gamma_1(2)$. Thus:
\begin{equation}
 [ SL_2(\Z):\Gamma_1(n)]=
\left\{ \begin{array}{cl}
n^2\prod_{p|n \  prime}\left(1-\frac{1}{p^2}\right) & \textrm{ if }\  n>2 \\
3 & \textrm{ if } \  n=2
\end{array} \right.
\end{equation}

\section{Quadratic growth rates and constants on marked tori}


\subsection{The 2 marked torus $\T^2_x$}

We begin with a general observation concerning the directions of periodic
families on the torus. These directions are exactly all the ``rational" directions and
this fact will not change how many points we will mark and wherever  
we will place them on the torus. The only thing that might change 
is the number of periodic families in a given rational direction, but in this case 
they are all of equal length. 
For the two marked torus $\T^2_x$ there are always two families in any 
rational direction if $x$ is non rational. There are saddle connections 
connecting $[0]$ und $[x]$ but which bound no periodic family because they 
are not in a rational direction. For rational markings $x$ this can never happen.

Collecting things for rational $x$: either there is a rational direction with two 
periodic families and the bounding saddle connections connecting the same marked 
point, or there is only one family but then the two bounding saddle connections
will connect the marked points $[0]$ and $[x]$.
Therefore counting the length $ N_{po}(\T^2_x,T)$ of cylinders of periodic trajectories   
for extremal markings $x$ on $\T^2_x$ is easy: 
\begin{equation}\label{ponr} 
\lim_{T \rightarrow \infty} \frac{N_{po}(\T^2_x,T)}{T^2}=
2\lim_{T \rightarrow \infty} \frac{N_{po}(\T^2,T)}{T^2}= \frac{6}{\pi} 
\quad \forall x \notin \T^2_{rat} 
\end{equation}
To treat the rational case we have to subtract the directions where 
only one family is from the doubled number of periodic orbits of the 
one marked torus. Let us denote the number of single cylinders on $\T^2_x$ 
with length smaller than $T$ by $N_{po,[x]}(\T^2_x,T)$. 
Thus we write: 
\[N_{po}(\T^2_x,T)=2N_{po}(\T^2,T)-N_{po,[x]}(\T^2_x,T)+o(T^2)\]
With the help of the discussion at the end of the last paragraph
one computes for a rational point $x=(\frac{p_1}{n},\frac{p_2}{n})$:
\begin{eqnarray}\label{po}
 \lim_{T \rightarrow \infty} \frac{N_{po}(\T^2_x,T)}{T^2}&=&
\lim_{T \rightarrow \infty} \frac{1}{T^2}\left(2N({\Z^2}^{V},T)
-N(G^{V}_{I^n},T)\right)=\nonumber\\
&=&\frac{6}{\pi}\left( 1-\frac{1}{2n}\prod_{p|n \ prime}
\left(1+\frac{1}{p}\right)^{-1}\right)
\end{eqnarray}
To count saddle connections is less easy, but again using the last paragraph
(and the discussion before equation (\ref{kxn}) we can write down the associated  
growth rates as follows:  
\begin{eqnarray*}
& &\lim_{T \rightarrow \infty}\frac{N_{sc}(\T^2_x,T)}{T^2}=\\
& &\lim_{T \rightarrow\infty}\frac{1}{T^2}\left(2N({\Z^{2}}^{V},T)
-2N(G^{V}_{I^n},T)\right)+\\
&+&\lim_{T \rightarrow\infty}\frac{1}{T^2}\sum_{i\in I^n}\left(\frac{1}{i^2}+
\frac{1}{(n-i)^2}\right)
\frac{n^2}{|I^n|}N(G^{V}_{I^n},T) = \cdots
\end{eqnarray*}
Here (\ref{vffinal}) and Proposition \ref{invariante} is used. The
last term needs some explanation. We have seen that $N(G^{V}_{I^n},T)$ already counts the length of the 
simple periodic cylinders. It counts as well the length of their boundaries
that is the sum of the length of the two bounding saddle connections.
By Proposition \ref{invariante} there are $\frac{1}{2}\varphi(n)$
different families of orbits of saddle connections labeled by their length ratios.
Equation (\ref{kxn}) shows that these orbits under the Veech group all have the 
same growth rate, namely $\frac{2}{\varphi(n)} N(G^{V}_{I^n},T)$.
The relation of the length of a saddle connection bounding a family to the 
length of the family is more or less the invariant $I$. Anyway it is one 
of the numbers $\frac{i}{n}$ or $\frac{n-i}{n}$ where $i \in I^n$. 
With the help of equation (\ref{scaling}) one finds the above expression. 
Continuing the evaluation: 
\begin{eqnarray}\label{rat}
\cdots &=& \frac{6}{\pi}-\frac{6}{n\pi}\prod_{p|n \ prime}
\left(1+\frac{1}{p}\right)^{-1}+
\frac{6}{\pi}\prod_{p|n \ prime}
\left(1-\frac{1}{p^2}\right)^{-1}\left(\sum_{i\in I^n}
\frac{1}{i^2}\right)=\nonumber \\
&=&\frac{6}{\pi}\left(1 + \prod_{p|n \ prime}\left(1+\frac{1}{p}\right)^{-1}
\left(\prod_{p|n \ prime}\left(1-\frac{1}{p}\right)^{-1}\sum_{i\in I^n}
\frac{1}{i^2}-\frac{1}{n}\right)\right)=\nonumber\\
&=&\frac{6}{\pi}\left(1 + \prod_{p|n \ prime}\left(1-\frac{1}{p^2}\right)^{-1}
\sum_{i\in I^n}\left(\frac{1}{i^2}-\frac{1}{n^2}\right)\right)
\end{eqnarray} The last equality follows with the help of equation (\ref{in}). 
For non rational markings one obtains:
\begin{eqnarray} \label{snrat} 
N_{sc}(\T^2_x,T)&=&2N({\Z^{2}}^{V},T) + N(x\Z^2,T)+o(T^2)=\nonumber \\
&=& \frac{6}{\pi}T^2+{\pi}T^2 +o(T^2).
\end{eqnarray}
Summarizing:
\begin{satz}The limits
\[ \lim_{T \rightarrow \infty} \frac{N_{po/sc}(\T^2_x,T)}{T^2}\] 
exist for all $x$, moreover there are estimates
\[\lim_{T \rightarrow \infty} 
\frac{N_{sc}(\T^2_x,T)}{T^2}\leq \lim_{T \rightarrow \infty} 
\frac{N_{sc}(\T^2_y,T)}{T^2}=\frac{6}{\pi}+\pi\] 
and
\[\lim_{T \rightarrow \infty} 
\frac{N_{po}(\T^2_x,T)}{T^2}\leq \lim_{T \rightarrow \infty} 
\frac{N_{po}(\T^2_y,T)}{T^2}=\frac{6}{\pi}\]
for  $y$ non rational. The inequalities  are strict if $x$ is rational.
 Moreover for each sequence of rational points $\{x_n\}^{\infty}_{n=0}$ converging to a  
 non rational point $x$:
\[\lim_{n\rightarrow \infty}\left(\lim_{T \rightarrow \infty} 
\frac{N_{sc/po}(\T^2_{x_n},T)}{T^2}\right)= \lim_{T \rightarrow \infty} 
\frac{N_{sc/po}(\T^2_x,T)}{T^2}\]
This is the continuity of the quadratic constants at non rational points. 
\end{satz}
\begin{proof}
The existence of the limits is proved in the equations (\ref{ponr}), (\ref{po}), (\ref{rat}) 
and (\ref{snrat}), even more the computations give explicitly the values of the limits. 
In the case of saddle connections it 
is indeed not easy to see the continuity of the growth rate function 
from the expression (\ref{rat}) directly (in fact before the author found the esimate
below he does not believe in the stated formula). 
By comparing equation (\ref{snrat}) (the limiting 
constant) to the first line of (\ref{rat}) one has to show:
\[\lim_{n \rightarrow \infty}\prod_{p|n \ prime}
\left(1-\frac{1}{p^2}\right)^{-1}\left(\sum_{i\in I^n}
\frac{1}{i^2}\right)= \frac{\pi^2}{6}\]
By using the equation (\ref{sumform}), which states
\[\sum_{i \in I^n}\sum^{\infty}_{l=0}\frac{1}{(nl+i)^2}=\frac{\pi^2}{6}\prod_{p|n \ prime}
\left(1-\frac{1}{p^2}\right)\]
and the estimates 
\[\sum_{i \in I^n}\left(\frac{1}{i^2}+\frac{\pi^2}{6}\frac{1}{4n^2}\right)<
\sum_{i \in I^n}\sum^{\infty}_{l=0}\frac{1}{(nl+i)^2}<
\sum_{i \in I^n}\left(\frac{1}{i^2}+\frac{\pi^2}{6}\frac{1}{n^2}\right)\]
we have:
\begin{eqnarray}\frac{\pi^2}{6}\left(1-\frac{1}{n}\prod_{p|n \ prime}
\left(1+\frac{1}{p}\right)^{-1}\right)<\sum_{i \in I^n}\frac{1}{i^2}\prod_{p|n \ prime}
\left(1-\frac{1}{p^2}\right)^{-1}< \nonumber \\
<\frac{\pi^2}{6}\left(1-\frac{1}{4n}\prod_{p|n \ prime}
\left(1+\frac{1}{p}\right)^{-1}\right).
\end{eqnarray}
Clearly if $n\rightarrow \infty$ the right and left hand side of the inequalities converge to
$\frac{\pi^2}{6}$.
The continuity in the case of periodic orbits follows immediately from
equation (\ref{po}).
\end{proof}
%

\subsection{The general case}


We generalize the results of the last section to arbitrary many markings, 
with the exception that we do not try to find explicit formulas for markings 
which are not in general position.
\begin{satz}\label{tlimes}
\begin{itemize}
\item[1.] [continuity] Let $\T^2_{[x_0,...,x_{n-1}]}$ be a n marked translation torus. 
Then the limits 
\[\lim_{T \rightarrow \infty}\frac{N_{sc/po}(\T^2_{[x_0,...,x_{n-1}]},T)}{T^2}\] 
exist.
\item[2.] [ $c^{max}_{sc/po}$ on markings in general position]
The set $\T^{2n}_{gp(sc/po)}$ has full Lebesque measure on $(\T^2)^{n}$ 
and the inequalities
\[\lim_{T \rightarrow \infty}\frac{N_{sc}(\T^2_{[x_0,...,x_{n-1}]},T)}{T^2}
\leq \lim_{T \rightarrow \infty}\frac{N_{sc}(\T^2_{[y_0,...,y_{n-1}]},T)}{T^2}=
\frac{n(n+1)}{2}\pi+\frac{3n}{\pi}\] and
\[\lim_{T \rightarrow \infty}\frac{N_{po}(\T^2_{[x_0,...,x_{n-1}]},T)}{T^2}
\leq \lim_{T \rightarrow \infty}\frac{N_{po}(\T^2_{[y_0,...,y_{n-1}]},T)}{T^2}=
 \frac{3n}{\pi}\] hold for all $(y_0,...,y_{n-1})\in \T^{2n}_{gp(sc/po)}$. 
 In both cases ($sc$'s and $po$'s) the inequality is strict, if 
 $(x_0,...,x_{n})\in \T^{2(n+1)}_{rat(sc/po)}$
\item[3.] [continuity at points in general position] \\
If $\ \{(x^i_0,...,x^i_{n})\}^{\infty}_{i=0} \quad  (x^i_0,...,x^i_{n})
\in \T^{2n}$ is a sequence with \\
$\lim_{i\rightarrow \infty}(x^i_0,...,x^i_{n})$ $=(y_0,...,y_{n-1})\in 
\T^{2n}_{gp(sc/po)}$
then 
\begin{equation}
\lim_{i\rightarrow \infty}\left(\lim_{T \rightarrow \infty}
\frac{N_{sc/po}(\T^2_{[x^i_0,...,x^i_{n}]},T)}{T^2}\right)= 
\lim_{T \rightarrow \infty}\frac{N_{sc/po}(\T^2_{[y_0,...,y_{n}]},T)}{T^2}
\end{equation}
\end{itemize}
\end{satz}
\begin{proof}
We prove the statements by induction over the number $n$ of marked points. 
The idea is to write down, modulo lower order terms, the numbers of saddle connections (or po's) for a given marking 
as sums over the numbers for an marking with one marked point less and  
the terms occuring when one marks the forgotten point. The results from the first 
section will imply the statements.
To start for one point there is nothing to show, moreover we have already seen
the proof for $n=2$. So we assume the statements are true for all n markings
$(x_0,...,x_{n-1})$ of the torus. We add another marking $x_n$. 
If the resulting marking is in general position, $(x_0,...,x_n) \in \T^{2(n+1)}_{gp(sc)}$
(for saddle connections) we obtain for saddle connections
\[N_{sc}(\T^2_{[x_0,...,x_{n}]},T)=N_{sc}(\T^2_{[x_0,...,x_{n-1}]},T)+
\sum^{n}_{i=0}N(SC_{(x_i,x_n)},T)\]
furthermore by the definition of general position markings 
\[\lim_{T \rightarrow \infty} \frac{N(SC_{(x_i,x_n)},T)}{T^2}=\pi\]
for all $i \in \{0,1,...,n-1\}$. For the case  $i=n$: saddle connections 
starting and ending in $x_n$ are bounding cylinders of periodic trajectories.
We have to take care of multiplicities when counting these
\[\lim_{T \rightarrow \infty} \frac{N(SC_{(x_n,x_n)},T)}{T^2}=\frac{3}{\pi}.\]
Thus
\begin{eqnarray}
\lim_{T \rightarrow \infty}\frac{N_{sc}(\T^2_{[x_0,...,x_{n}]},T)}{T^2}&=&
\lim_{T \rightarrow \infty}\frac{N_{sc}(\T^2_{[x_0,...,x_{n-1}]},T)}{T^2}+
n\pi+\frac{3}{\pi}=\nonumber \\
&=& \frac{n(n+1)}{2}\pi+\frac{3n}{\pi}. \nonumber
\end{eqnarray}
We assume now that the marking is in general position with respect to periodic orbits.
Since a saddle connection starting and ending in $x_n$ bounds a cylinder of periodic
trajectories, one has to add in each rational direction 
a new periodic family. Thus the periodic orbit growth rates for markings in general position are:
\[\lim_{T \rightarrow \infty}\frac{N_{po}(\T^2_{[x_0,...,x_{n}]},T)}{T^2}=
\lim_{T \rightarrow \infty}\frac{N_{po}(\T^2_{[x_0,...,x_{n-1}]},T)}{T^2}+
\frac{3}{\pi}=\frac{3n}{\pi}.\]
If the marking $(x_0,..,x_n)$ is not in general position, it could be divided in 
classes of relatively rational points (this is of course meaningful in any case, po's or sc's).  
The set of these classes is denoted by $\mathcal{R}_{[x_0,...,x_n]}$.
For each class the Veech theory holds
(see \cite{ve89,ve92, ve98}) and guaranties that all the limits exist.
For the inequalities between the limits one observes: if the marked points
are rational there exists always a direction in which there are strictely less than 
the maximal number of $n$ cylinders of closed trajectories. 
The orbit of this direction under the operation of the
Veech group leads to a positive asymptotic growth constant (Proposition 6.1 in \cite{gj97}). 
But a lower (than maximal possible) number of closed cylinders in a given direction 
causes ``multiplicities" (on has to count only visible points). 
Hence the limit growth rates are in both cases (po's and sc's) smaller than the ones for  markings in general position. 
This shows: if ther is one class in $\mathcal{R}_{[x_0,...,x_n]}$ containing more than one point,
then the growth rates have to be smaller than for markings in general position.
To compute the growth rates in the cases of periodic orbits, we observe that
each class $m\in \mathcal{R}_{[x_0,...,x_n]}$ defines a marked Veech torus $\T^2_m$, 
by marking $\T^2$ only with the points of that class. Thus:
\begin{equation}
N_{po}(\T^2_{[x_1,...,x_{n}]},T))=\sum_{m \in \mathcal{R}_{[x_0,...,x_n]} }
N_{po}(\T^2_{m},T) +o(T^2).
\end{equation}
Here is used that by definition between points of different classes there are never
saddle connections with rational slope (or equivalently: with bound a closed trajectory).

All remains is to count the saddle connections between points of different classes. 
To do this let $x_i$ be a  marked point not rationally equivalent to $x_n$ and 
$SC^{(x_0,...,x_{n-1})}_{(x_i,x_j)}$ the set of saddle connections connecting 
$x_i$ and $x_j$ on the n marked torus $\T^2_{[x_1,...,x_{n-1}]}$.
What we have to count are sets of saddle connections from a given set of saddle connections 
between two marked points (for example the set $SC^{(x_0,...,x_{n-1})}_{(x_i,x_j)}$) 
which are parallel and shorter to another set of saddle connections (for example $SC_{(x_i,x_n)}$).
However by Lemma (\ref{schnigi}) and Proposition (\ref{countnrat}) 
sets of this kind are always lattices modulo a mistake of order $O(T)$, which does not count 
in the quadratic limit.\\  
There are two different cases: First, there is a nonempty subset in 
$SC_{(x_i,x_n)}$  consisting of saddle connections 
parallel to saddle connections in the set $SC^{(x_0,...,x_{n-1})}_{(x_i,x_j)}$ and shorter in length. 
The set of all indices $j$ of such $x_j$ are denoted with $I$. 
Second,  a nonempty subset $SC_{(x_i,x_n)}$ is parallel to one 
from $SC^{(x_0,...,x_{n-1})}_{(x_i,x_j)}$, but the first are longer. 
As above denote the set of all indices of such $x_j$ with $K$. 
Moreover the sets of saddle are (modulo sets of lower order growth rates) 
given by some lattices $G^{(x_i,x_n)}_{(x_i,x_j)}$ or $G^{(x_i,x_j)}_{(x_i,x_n)}$ respectively. 
We can write:  
\begin{equation}
N(SC^{(x_0,...,x_n)}_{(x_i,x_n)},T)=N(G_{(x_i,x_n)},T)-\sum_{x_j\in I}
N(G^{(x_i,x_j)}_{(x_i,x_n)},T) +o(T^2).
\end{equation}
Here $G_{(x_i,x_n)}$ denotes the lattices of saddle connections between $x_i$ and $x_n$ on the 
two marked tori $\T^2_{[x_i,x_n]}$. For all $j\in K$
\begin{equation}
N(SC^{(x_0,...,x_n)}_{(x_i,x_j)},T)=N(SC^{(x_0,...,x_{n-1})}_{(x_i,x_j)},T)
-N(G^{(x_i,x_n)}_{(x_i,x_j)},T)+o(T^2).
\end{equation}
Modulo terms of strictly lower order size the above expressions are all lattices, 
hence they have quadratic growth rate limits.
After these steps for all markings $[x_i]$ the induction is complete. 
From the two equalities it follows immediately:\\
if intersections of the above kind have a nonzero quadratic growth rate 
(so the marking is rational in the sense of our definition) than the associated growth rates are 
smaller than for a non rational marking. The continuity follows from the fact that by approximating 
a marking in general position all the denominators of all the length ratios are growing over all bounds.  
\end{proof}


\subsection{Branched coverings of marked tori}

Because of the simplicity of our methods we are not able to conclude anything about (families of) 
branched coverings of tori, without using the existence of the growth rate limits by the  
method explained in \cite{ems}. Even the sets where the growth rate
is the biggest is not seen as just the preimage of the non rational markings on
the torus. For closed cylinders of geodesics at least it has to be contained in this set. 
This is because the geometry and combinatorics of the covering might cause new sets 
where the growth rates are different.
Geometrically the reason for this is that the length spectrum  of the cylinders of periodic
trajectories is connected in a way to the saddle connections bounding them, which depends on the 
covering and where the singular points are.  
If one makes the simplifying assumption that all inverse images of the marked points
are itself marked (such coverings are called `balanced coverings' by some people), then 
from Theorem \ref{tlimes} one can conclude:
\begin{kor}
Let $\V \stackrel{\pi}{\rightarrow} \T^2_{[x_0,...,x_n]}$ a covering   
branched only over $(x_0,...,x_n)$ with the induced translation structure 
$\pi^{\ast}dz$ . Further denote by ${\V}_C$ the translation structure 
where each inverse image of a marked point is itself marked . 
Then, after rescaling the volume of $ \V$ (with respect to the induced flat 
metric) to one:
\begin{equation}
\lim_{T \rightarrow \infty}\frac{N_{sc}(\V_C,T)}{T^2}=
deg(\pi)^2 \lim_{T \rightarrow \infty}\frac{N_{sc}(\T^2_{[x_0,...,x_{n}]},T)}{T^2}
\end{equation}
($deg(\pi)$ is the degree of $\pi$). Moreover, if some of the markings are artificial
(i.e. not singular points of the induced metric) then the growth rates after ``unmark" 
any subset of these points will be smaller than the above one. Thus
\begin{equation}\label{verh}
\limsup_{T \rightarrow \infty}\frac{N_{sc}(\V,T)}{T^2}\leq
deg(\pi)^2 \left(\frac{n(n+1)}{2}\pi+\frac{3n}{\pi} \right).
\end{equation}
\end{kor}
Remark: Similarly one can write down the corresponding estimate for maximal periodic cylinders. \\
If one includes Ratners Theorem as done in \cite{ems} then
the $\limsup$ in equation (\ref{verh}) can be replaced by a limit.\\
\begin{proof}
The statement is a direct consequence of the fact that each saddle connection 
on the torus has exactly $deg(\pi)$ preimages on $\V_C$ and Theorem \ref{tlimes}.
The second part is clear, because by removing nonsingular markings, 
the directions in which there are saddle connections are not changed.  
But in any such direction the number of saddle connections could decrease and their length increase.
\end{proof}
Remark: With the results of Veech \cite{ve90} every translation surface can
be approximated by torus coverings in spaces of abelian differentials ${\cal A}(g,P)$. 
Moreover torus coverings are equally distributed with respect to the 
Liouville measure $\mu$ mentioned in the introduction.
Since periodic cylinders and saddle connections are locally stable with 
respect to deformations in this spaces we have some sort of control
on their numbers by the Corollary. But because of the increasing
number of artifical markings while approximating a general translation surface
or equivalently the increasing degree of the covering, the estimate above
is to weak to predict for example quadratic bounds of the growth rates of 
po's or sc's in general. On the other hand the equidistribution of torus coverings in ${\cal A}(g,P)$
is used to compute the $\mu$ volume of these spaces of abelian differentials see \cite{eo}. 
The knowledge of these volumes in turn is the main step to evaluate Siegel Veech formulas, compare 
\cite{ems, emz}.

\vspace*{2cm}
\begin{center} 
Authors present adress:
\vspace*{5mm}\\
Technische Universit\"at Berlin \\ Strasse des 17. Juni 136 \\  
Sekr. MA 7-1/ MA 7-2 \\ D-10623 Berlin \\ 
\normalsize{e-mail: \texttt{schmoll@math.tu-berlin.de} }

\end{center}

\end{document}